\documentclass[letter,reqno]{amsart}
\usepackage{amsmath,amsfonts,amsthm, amssymb}
\usepackage{color}
\usepackage{graphicx}

\theoremstyle{plain}

\newtheorem*{lemma}{Lemma}

\newtheorem*{corollary}{Corollary}

\theoremstyle{remark}
\newtheorem*{remark}{Remark}

\newcommand{\R}{\mathbb{R}}
\newcommand{\dd}{\mathrm{d}}

\usepackage{hyperref}

\begin{document}

\title[A replacement lemma for obtaining pointwise estimates]{A replacement lemma for obtaining pointwise estimates in phase transition models}
\author{Nicholas D.\ Alikakos}
\author{Giorgio Fusco}
\address{Department of Mathematics\\ University of Athens\\ Panepistemiopolis\\ 15784 Athens\\ Greece \and Institute for Applied and Computational Mathematics\\ Foundation of Research and Technology -- Hellas\\ 71110 Heraklion\\ Crete\\ Greece}
\email{\href{mailto:nalikako@math.uoa.gr}{\texttt{nalikako@math.uoa.gr}}} 
\address{Dipartimento di Matematica Pura ed Applicata\\ Universit\`a degli Studi dell'Aquila\\ Via Vetoio\\ 67010 Coppito\\ L'Aquila\\ Italy} \email{\href{mailto:fusco@univaq.it}{\texttt{fusco@univaq.it}}}

\begin{abstract}
We establish a replacement lemma for a variational problem, which is not based on a local argument. We then apply it to a phase transition problem and obtain pointwise estimates.
\end{abstract}

\maketitle

\section{Introduction}
We consider the elliptic system
\begin{equation}\label{system}
\Delta u - W_u(u) = 0, \text{ for } u: \Omega \subset \R^n \to \R^m,
\end{equation}
where $W : \R^m \to \R$ a nonnegative $C^1$ potential possessing several minima and $W_u (u) := ( \partial W/ \partial u_1, \ldots, \partial W / \partial u_n )^{\top}$. The system \eqref{system} is variational with associated functional
\begin{equation}\label{action}
J_\Omega(u) = \int_\Omega \left(\frac{1}{2}|\nabla u|^2+W(u) \right) \dd x.
\end{equation}
In what follows, we take $\Omega$ to be a bounded, open, and connected set in $\R^n$, with Lipschitz boundary. We introduce the hypothesis

\medskip

\noindent {\bf (H)} {\em Let $\lambda \to W(a+ \lambda w)$, with $|w| = 1$, be a strictly increasing function on $[0, r_0)$. The vector $a$ is a global minimum of $W$ and $r_0$ is positive and fixed.}

\medskip
\noindent Note that (H) is a very weak nondegeneracy hypothesis that was introduced in \cite{4}.

\medskip

The main purpose of this note is to establish the following

\begin{lemma}
Let $\Omega$ be as above and let $A \subset \Omega$ be an open, Lipschitz set with $\partial A \neq \varnothing$. Moreover, suppose that 

\begin{enumerate}
\item $u(\cdot)\in W^{1,2}(\Omega)\cap C^1(\Omega)$,\medskip
\item $|u(x)-a| \leq r$ on $\partial A\cap\Omega$, for some $r$ with $2r\in(0, r_0)$,\medskip
\item there is an $x_0\in A$ such that $|u(x_0)-a|>r$.\medskip
\end{enumerate}

Then, there exists $\tilde{u}(\cdot)\in W^{1,2}(\Omega)$ such that
\begin{equation*}
\left\{
\begin{array}{ll}
\tilde{u}(x) = u(x), & \text{ in }  \Omega \setminus A, \medskip\\
|\tilde{u}(x)-a| \leq r, & \text{ in } A, \medskip\\
J_\Omega (\tilde{u}) < J_\Omega(u). 
\end{array}\right.
\end{equation*}

\end{lemma}

\begin{figure}
\begin{picture}(0,0)%
\includegraphics{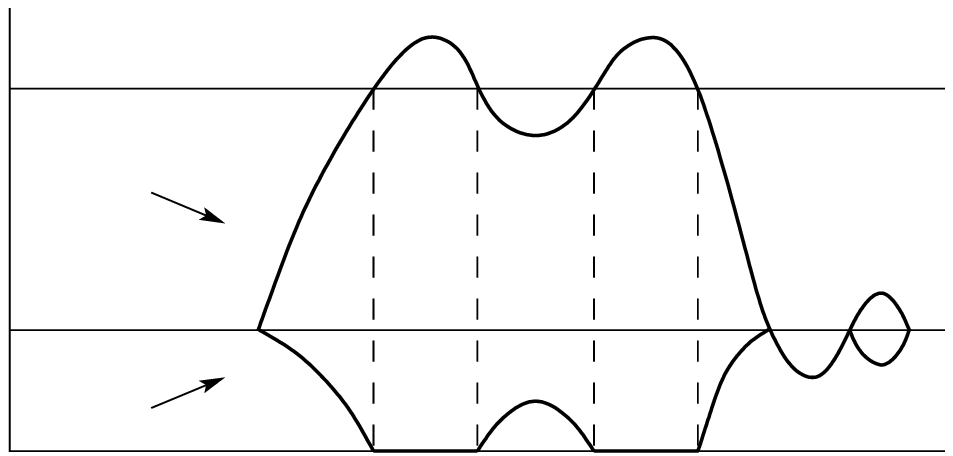}%
\end{picture}%
\setlength{\unitlength}{4144sp}%
%
\begin{picture}(4597,2059)(3403,-4020)
\put(3491,-3464){$r$}
\put(4191,-2815){$u$}
\put(4191,-3860){$\tilde{u}$}
\put(3463,-2366){$r_0$}
\end{picture}%
\caption{}
\label{figure1}
\end{figure}

We note that in the lemma, no {\em a priori} bound is imposed on the $\max_A |u(x) - a|$ and, thus, the lemma is not of local nature. Its meaning is that from the point of view of minimizing $J$ for a function that is in part close to the minimum value of $W$, independently of the structure of $W$, it is more efficient to remain close to the minimum throughout (see Figure \ref{figure1}).

We will illustrate the lemma above by establishing the following pointwise estimate.

\begin{corollary}
Let $n = m = 2$ and let $W$ have exactly one global minimum at $a =( \alpha, 0)$ on the right half-plane $\R^{2}_{+} = \{ u_1, u_2 \mid u_1 \geq 0 \}$, while $W>0$ in $\R^{2}_{+} \setminus \{a\}$. Consider the family of variational problems
\[
\min J_{\Omega_R^\mu}, \text{ where } \Omega_R^\mu = \{(x_1, x_2) \mid 0< x_1 <\mu R \text{ and } |x_2| \leq R\},
\]
with corresponding global minimizers $\{u_{R,\mu}\}$ and suppose that
\begin{enumerate}
\medskip \item[(i)] $u_{R,\mu}$ maps $\Omega^\mu_R$ in $\{(u_1, u_2) \mid u_1\geq 0\}$, \quad (positivity)
\medskip \item[(ii)] $J_{\Omega_R^\mu} (u_{R,\mu}) \leq CR$, where $C$ a universal constant,
\medskip \item[(iii)] $|u_{R,\mu}(x) - a| \leq r$ on $\{(\mu R, x_2) \mid |x_2| \leq R\}$ and $\partial
u_{R,\mu} / \partial n = 0$ on the remaining three sides of $\partial\Omega^\mu_R$.
\end{enumerate}

Then, there exist $R_0 > 0$, $\mu_0 > 0$, and $\eta_0 > 0$ such that
\[
|u_{R,\mu}(x)-a| \leq \frac{r}{2} ~\text{ in }~ \{(x_1, x_2) \in \Omega^\mu_R \mid \eta_0 R \leq x_1 \le\mu R\},
\]
for all $R \geq R_0$ and $\mu \geq \mu_0$.
\end{corollary}

The proof of the corollary is a two-dimensional measure-theoretic argument, where the kinetic and potential terms in the energy are estimated independently. It would be very interesting to extend this to higher dimensions. The one-dimensional version of the lemma above appeared in \cite{4}, and subsequently in \cite{5}, where an extension from balls to convex sets was given. For hypotheses (i) and (ii) see \cite{1}, \cite{3}.

\section{Proofs}

\begin{proof}[Proof of the Lemma.]
We utilize the polar representation
\begin{equation}\label{polar}
u(x) = a + | u(x ) -a | \frac{u(x) - a}{ | u(x) - a | } =: a + \rho(x) n(x)
\end{equation}
and note that
\[
|\nabla u(x)|^2 = |\nabla \rho(x)|^2+\rho^2(x) |\nabla u(x)|^2.
\]

\medskip

{\em Step 1.} We begin by settling the lemma under the additional hypothesis
\begin{equation}\label{additional}
\rho(x) \leq 2r < r_0, \text{ in } A.
\end{equation}
We choose $\varepsilon > 0$ so that
\begin{equation}\label{critical}
\rho(x)>r+\varepsilon, \text{ where }r+\frac{\varepsilon}{2} \text{ is not a critical value
of } \rho \text{ in } A.
\end{equation}
Therefore, the set
\[
\Gamma_\varepsilon = \partial C^1_\varepsilon \cap A, \text{ where } C_\varepsilon = \left\{ x \in
A \mid \rho(x) > r + \frac{\varepsilon}{2} \right\},
\]
is a $C^1$ manifold in $A$.

Now, define $\tilde{u}_\varepsilon$ as follows.
\begin{equation}\label{u-tilde}
\left\{\begin{array}{ll}
\tilde{u}_\varepsilon (x) = u(x), & \text{in } A \setminus C_\varepsilon,\smallskip\\
\tilde{u}_\varepsilon (x) = a + \left( r + \dfrac{\varepsilon}{2} \right) n(x), & \text{in } C_\varepsilon,\smallskip\\
\tilde{u}_\varepsilon (x) = u(x), & \text{in } \Omega \setminus A.
\end{array}\right.
\end{equation}
Notice that $\tilde{u}_\varepsilon$ is continuous on $\Gamma_\varepsilon$. There also holds
\[
|\nabla\tilde{u}_\varepsilon (x) |^2 = \left( r + \frac{\varepsilon}{2} \right)^2 |\nabla n(x)|^2 \leq \rho^2(x) |\nabla n(x)|^2 \leq |\nabla u(x)|^2
\]
in $C_\varepsilon$. It follows that $\tilde{u}_\varepsilon \in W^{1,2}(\Omega)$ and,
moreover,
\begin{equation}\label{energy-estimate}
\int_\Omega |\nabla u|^2 \,\dd x \geq \int_\Omega |\nabla\tilde{u}_\varepsilon|^2 \,\dd x.
\end{equation}
Hence, $\tilde{u}_\varepsilon \rightharpoonup \tilde{u}$ in $W^{1,2}$ as $\varepsilon \to 0$, and by weak lower semi-continuity,
\begin{equation}\label{energy-estimate-limit}
\int_\Omega |\nabla u|^2 \,\dd x \geq \int_\Omega |\nabla\tilde{u}|^2 \,\dd x.
\end{equation}
Clearly
\[
\left\{\begin{array}{ll}
\tilde{u}(x) = a + r n(x), & \text{in } C_0 = \{x \in A \mid \rho(x) > r \}, \medskip\\
\tilde{u} = u(x), & \text{in } \Omega \setminus C_0.
\end{array}\right.
\]
Finally,
\[
\int_A W(u(x)) \,\dd x = \int_{A \setminus C_0} W( a + \rho(x) n(x) ) \,\dd x + \int_{C_0} W( a + \rho(x) n(x) ) \,\dd x.
\]
By (H), (iii), and the hypothesis $A^+ = \{ x \in A \mid \rho(x) > 2r \} = \varnothing$,
\[
\int_{C_0} W( a + \rho(x) n(x) ) \,\dd x > \int_{C_0} W ( a + r n(x) ) \,\dd x.
\]
Therefore,
\begin{equation} \label{w-estimate}
\int_\Omega W(u) \,\dd x > \int_\Omega W(\tilde{u}) \,\dd x, 
\end{equation}
and so, $J_\Omega (u) > J_\Omega (\tilde{u})$.

Also by (\ref{u-tilde}),
\[
\left\{\begin{array}{ll}
\tilde{u}(x) = u(x), & \text{in } A \setminus C_0, \medskip\\
\tilde{u}(x) = a + r n(x), & \text{in } C, \medskip\\
\tilde{u}(x) = u(x), & \text{in } \Omega \setminus A, \\
\end{array}\right.
\]
thus, the lemma is established under hypothesis \eqref{additional}.

\medskip

{\em Step 2.} We may therefore assume that
\begin{equation}\label{A+}
|A^+| > 0.
\end{equation}
We first assume that $r$ is not a critical value of $\rho$ in $A$ and later we remove this assumption.

Define the Lipschitz function
\begin{equation}\label{alpha-tau}
\alpha (\tau) = \left\{\begin{array}{ll}
 1, & \text{for } \tau \leq r, \medskip\\
  \dfrac{2r-\tau}{r}, & \text{for } r \leq \tau \leq 2r, \medskip\\
  0, & \text{for } r \geq 2r, \\
\end{array}\right.
\end{equation}
and recall that compositions of Lipschitz functions with $W^{1,2}$ functions render $W^{1,2}$ functions.

Set
\begin{equation}
\left\{\begin{array}{ll}
w(x) = u(x), & \text{in } A \setminus C_0 \medskip\\
w(x) = a + r \alpha (\rho(x)) n(x), & \text{in }  C \medskip\\
w(x) = u(x), & \text{in } \Omega \setminus A. \\
\end{array}\right.  \label{eq12}
\end{equation}
Note that $W$ is continuous on $\partial C$ ($C^1$ manifold) and so $w$ is in $W^{1,2}(\Omega)$.

In $\{ x\in A \mid r \leq \rho(x) \leq 2r \}$ there holds
\begin{align}\label{estimate1}
|\nabla w(x)|^2 &= |\nabla\rho(x)|^2 + r^2 \alpha^2 |\nabla n(x)|^2 \nonumber \\
&\leq |\nabla\rho(x)|^2 + r^2 |\nabla n(x)|^2 \quad \text{(since
$\alpha \leq 1$)} \nonumber\\
&\leq |\nabla\rho(x)|^2 + \rho^2 |\nabla u(x)|^2 \nonumber \\
&= |\nabla u(x)|^2.
\end{align}

Also $\nabla w=0$ in $A^+$ and $\nabla w=\nabla u$ in the rest of $A$. It follows that
\begin{equation}\label{energy-uw}
\int_\Omega |\nabla u|^2 \,\dd x \geq \int_\Omega |\nabla w|^2 \,\dd x.
\end{equation}

In $\{x\in A \mid r\le\rho(x)\le2r\}$ there holds
\begin{align}\label{estimate2}
W(w(x)) &= W(a + r \alpha(\rho(x)) n(x)) \nonumber \\
&\leq W(a + r n(x)) \nonumber \\
&\leq W(a + \rho(x) n(x)) \nonumber \\
&= W(u(x)),
\end{align}
while
\[
W(w(x)) = 0 < W(u(x)), \text{ in } A^+,
\]
since $a$ is a global minimum. 

Now, since $|A^+|>0$, we obtain
\begin{equation}\label{uw-estimate}
\int_\Omega W(u(x)) \,\dd x > \int_\Omega W(w(x)) \,\dd x.
\end{equation}
We also note that
\[
| w(x) - a | \leq r, \text{ in } A.
\]

Thus, the lemma is established in this case as well.

\medskip

{\em Step 3.} Finally, suppose that $r$ is a critical value of $\rho$ in $A$. We can choose a decreasing and noncritical sequence $r_n \to r$. Then, the hypotheses (i), (ii), (iii) of the lemma are satisfied with $r=r_n$ and, thus, we obtain a sequence $\{\tilde{u}_n\}$ with the following properties:
\[
\left\{\begin{array}{ll}
\tilde{u}_n(x) = u(x), & \text{in } \Omega \setminus A, \medskip\\
|\tilde{u}_n(x) - a| \leq r_n, & \text{in } A, \medskip\\
J_\Omega(\tilde{u}_n) < J_\Omega(u).\\
\end{array}\right.
\]
Moreover, by construction,
\[
\int_\Omega |\nabla u|^2 \,\dd x \geq \int_\Omega |\nabla\tilde{u}_n|^2 \,\dd x.
\]
Hence, by taking possibly a subsequence, there holds $\tilde{u}_n \rightharpoonup \tilde{u}$ in $W^{1,2}$ as $n \to \infty$ and thus,
\[
\int_\Omega |\nabla u|^2 \,\dd x \geq \int_\Omega |\nabla\tilde{u}|^2 \,\dd x.
\]
By the compactness of the embedding $W^{1,2}_{\text{loc}} \hookrightarrow L^2_{\text{loc}}$ and from
\[
W(\tilde{u}_n(x))\le W(u(x)), \text{ in } \Omega,
\]
we obtain
\[
W(\tilde{u}(x))\le W(u(x)), \text{ a.e.\ in } \Omega.
\]
However,
\[
\int_{A^+} W(u) \,\dd x > \int_{A^+} W(\tilde{u}) \,\dd x,
\]
thus, it follows that
\begin{equation*}
\left\{
\begin{array}{ll}
\tilde{u}(x) = u(x), & \text{ in }  \Omega \setminus A, \medskip\\
|\tilde{u}(x)-a| \leq r, & \text{ in } A, \medskip\\
J_\Omega (\tilde{u}) < J_\Omega(u). 
\end{array}\right.
\end{equation*}

The proof of the lemma is complete.
\end{proof}

We continue with the

\begin{proof}[Proof of the Corollary.]
In what follows, we write $u$ for $u_{R,\mu}$, $\rho$ for $\rho_{R,\mu}$ etc. Consider the sets $j_R \subset i_R \subset \R$, with 
\[
i_R := \left\{ x_1 \in (0,\eta R) ~\Big|~ \text{there exists } x_2 \in (0,R) \text{ with }
\rho(x_1, x_2) \geq \frac{r}{2} \right\}
\]
and
\[
j_R := \left\{ x_1 \in i_R ~\Big|~ \text{there exists } x_2 \in (0,R) \text{ with }
\rho(x_1, x_2) \geq \frac{r}{4} \right\}
\]

Then, the positivity property (i) implies the lower bound
\begin{eqnarray}\label{lowerbound}
R w_0 | i_R \setminus j_R | \leq \int_{0}^{R} \int_{i_R \setminus j_R} W(u) \,\dd x_1 \dd x_2,
\end{eqnarray}
where $w_0 := \min_{|u-a|>r/4} W(u) > 0$.

From the definition of $j_R$, we conclude that for $x_1 \in j_R$ there is an interval $L_{x_1} = (a_{x_1}, b_{x_2})$ of $x_2$ values such that
\[
\frac{r}{4} = \rho(x_1, a_{x_1}) \leq \rho(x_1, x_2) \leq \rho(x_1, b_{x_2}) = \frac{r}{2}, \text{ for all } x_2 \in L_{x_1}.
\]
It follows that
\begin{equation} \label{w-w-estimate}
\int_{L_{x_1}} W(u(x_1,\tau)) \,\dd \tau \geq w_0 |L_{x_1}|, \text{ for all } x_1\in j_R.
\end{equation}
Moreover, we have
\begin{align}\label{r-L-estimate}
\frac{r}{4} &\leq \int_{L_{x_1}} \left |\frac{\partial \rho}{\partial x_2} (x_1,\tau) \right| \,\dd \tau 
\leq \left( |L_{x_1}| \int_{L_{x_1}} \left| \frac{\partial\rho} {\partial x_2}(x_1,\tau) \right|^2 \,\dd \tau \right)^{1/2} \nonumber\\
&\leq \left( |L_{x_1}| \int_{L_{x_1}}| \nabla u(x_1,\tau)|^2 \,\dd \tau \right)^{1/2}.
\end{align}
From (\ref{w-w-estimate}) and (\ref{r-L-estimate}) we have
\[
\frac{1}{32} \frac{1}{|L_{x_1}|} r^2 + w_0 |L_{x_1}| \leq \int_{L_{x_1}} \frac{1}{2}
|\nabla u(x_1,\tau)|^2 \,\dd\tau + \int_{L_{x_1}} W(u(x_1,\tau)) \,\dd\tau,
\]
thus,
\begin{equation}\label{r-sqrt}
\frac{r\sqrt{w_0}}{2\sqrt{2}} \leq \int_{L_{x_1}} \frac{1}{2} |\nabla u (x_1,\tau)|^2 \,\dd \tau + \int_{L_{x_1}} W(u(x_1,\tau)) \,\dd \tau.
\end{equation}

Concluding,
\begin{align}\label{conclusion}
CR &\overset{\text{(ii)}}{\geq} 2 \int_{\Omega_{R,\mu}} \left(\frac{1}{2} |\nabla u|^2 + W(u) \right) \dd x \geq \int^{R}_{0} \int_{i_R} \left(\frac{1}{2} |\nabla u|^2 + W(u) \right) \dd x_1 \dd x_2 \nonumber \\
&= \int^{R}_{0} \int_{i_R \setminus j_R} \left( \frac{1}{2} |\nabla u|^2 + W(u) \right) \dd x_1 \dd x_2 + \int^{R}_{0} \int_{j_R} \left( \frac{1}{2} |\nabla u|^2 + W(u) \right) \dd x_1 \dd x_2 \nonumber \\
&\geq R w_0 |i_R \setminus j_R| + \frac{r \sqrt{w_0}}{2\sqrt{2}} |j_R|,
\end{align}
where the last inequality follows from \eqref{lowerbound}, \eqref{r-sqrt}. Hence,
\begin{align}
CR &\geq A |j_R| + B (|i_R| - |i_R|)R, \text{ for } A:= r \sqrt{w_0} / 2 \sqrt{2},\ B := w_0, \nonumber \\
&\geq \min \{ A, BR \} |i_R| \nonumber \\
&\geq A |j_R|, \text{ if } R \geq r / 2 \sqrt{2w_0} =: R_0.
\end{align}

Consequently, if we take $R$ large, we obtain that 
\[
|i_R| \leq \frac{2\sqrt{2}CR}{r \sqrt{w_0}} =: \eta_0 R.
\]
If we take $\eta > \eta_0$ and fix it, then $|i_R| < \eta R$ and therefore there is an $\bar{x}_1 \in (0, \eta R)$, which does not belong to $i_R$, and such that
\begin{equation}\label{x1-bar}
\rho(\bar{x}_1, x_2) < \frac{r}{2}, \text{ for all } x_2 \in(0,R).
\end{equation}

Applying now the lemma for the choice $A = \{ (x_1, x_2) \mid \bar{x}_1 \leq x_1 \leq \mu R, |x_2| < R \}$, we conclude that $\rho \leq r/2$ in $A$, thus, $\rho < r$ on the line $x_1 = \eta R$.
\end{proof}

\begin{remark}
The intuition behind hypothesis (ii) is that if $u_{R,\mu}$ is bounded away from $a$ on a large set, then
\[
\int_{\Omega^\mu_R} W(u_{R,\mu}(x))\,\dd x \geq C R^2,
\]
therefore, by (ii) this cannot happen.

The {\em a priori} bound (ii) is related to the fact that \eqref{action} is linked to a perimeter functional (see \cite{2}). In general dimensions, the appropriate {\em a priori} estimate is $J_{\Omega_R} (u) \leq CR^{n-1}$.
\end{remark}


\nocite{*}
\bibliographystyle{plain}

\begin{thebibliography}{9}

\bibitem{1} S.\ Alama, L.\ Bronsard, and C.\ Gui. Stationary layered solutions in $\R^2$ for an Allen--Cahn system with multiple well potential. {\em Calc.\ Var.} {\bf 5} No. 4 (1997), pp.\ 359--390.

\bibitem{2} G.\ Alberti. Variational models for phase transitions, an approach via $\Gamma$-convergence. In {\em Calculus of variations and partial differential equations}, L.\ Ambrosio and N.\ Dancer, edited by G.\ Buttazzo, A.\ Marino, and M.\ K.\ V.\ Murthy. Springer, 2000.

\bibitem{4} N.\ D.\ Alikakos and G.\ Fusco. On the connection problem for potentials with several global minima. {\em Indiana Univ.\ Math.\ J.} {\bf 57} No.\ 4 (2008), pp.\ 1871--1906.

\bibitem{3} N.\ D.\ Alikakos and G.\ Fusco. On an elliptic system with symmetric potential possessing two global minima. Preprint. arXiv:0810.5009.

\bibitem{5} N.\ D.\ Alikakos and N.\ I.\ Katzourakis. Heteroclinic travelling waves of gradient diffusion systems. To appear in {\em Trans.\ Amer.\ Math.\ Soc.}

\end{thebibliography}

\end{document}